\documentclass[14pt]{extarticle}
\usepackage[margin=1in]{geometry}
\usepackage{dsfont}
\usepackage{amsfonts}
\newtheorem{theorem}{Theorem}[section]
\newtheorem{pro}{Proposition}[section]
\newtheorem{lemma}{Lemma}[section]

\newtheorem{remark}{Remark}[section]
\newtheorem{cor}{Corollary}[section]

\newcommand{\proof}[1]{\noindent{\it\bf Proof:#1\ }}
\newcommand{\QED}{\hfill$\Box$\medskip}

\begin{document}
\title{On   Nonlinear Part of Filled-Section in Splicing }
\author{
	\\
	Gang Liu\\Department of Mathematics\\UCLA }
\date{August,  2018}
\maketitle

\section{Introduction}
The purpose of this paper is to define the nonlinear part of the filled section $\Psi_N=(\Psi_{N, -},\Psi_{N, +} )$ denoted by $
N=(N_-, N_+) $ within the framework of the usual analysis of Banach manifolds rather than in the setting of Sc-analysis of polyfold theory in \cite{1, 2}. The main difficulty  for this is that the simple choice made in \cite{ 2} for $\Phi_-$ using a linear operator leads to a filled-section with loss of differentiability so that the theory for filled-sections has to be formulated in Sc-analysis.

To overcome this difficulty, recall that
for each fixed gluing parameter $a=(R, \theta) \in [R_0, \infty)\times S^1,$ 
$N_{\pm} ^a:L_{k,\delta}^p(C_-, E))\times L_{k,\delta}^p(C_+, E)\rightarrow L_{k-1,\delta}^p(C_-, E))\times L_{k-1,\delta}^p(C_+, E)$ is obtained from 
$\Phi_{\pm,N}^a:L_{k,\delta}^p(S^a_{\pm}, E)\rightarrow L_{k-1,\delta}^p(S^a_\pm, E)$
by the conjugation by the total gluing  map $T^a: C_-\cup C_+\rightarrow S^a$.  Here $C_{\pm}\simeq (0,\pm \infty)\times S^1$,  $E={\bf C}^n$ and $L_{k, \delta}^p$-maps  here    are the  $L_k^p$-maps  that decay exponentially along the ends of the half cylinders
$C_{\pm}$  with the decay rate $0<\delta<1.$

Since  $\Phi_+$ is required to be the usual ${\overline {\partial }}_J$-operator in the Gromov-Witten theory,  the nonlinear part $\Phi^a_{+, N}$ then is given by $\Phi^a_{+, N}(v_+ )=J(v_+)\partial_s v_+$, where $v_+:S^a_+\rightarrow {\bf C}^n$ is a $L_{k,\delta}^p$-map with the domain $S^a_+\simeq [-R,R]\times S^1$  with   cylindrical coordinate $(t,s)$.  To get desired $\Phi^a_N$ without loss of loss of differentiability, the key observation is that the choice of  $\Phi^a_{-, N}:L_{k,\delta}^p(S^a_-, E)) \rightarrow L_{k-1,\delta}^p(S^a_-, E))$ with $S^a_-\simeq {\bf R}^1\times S^1$ should be $\Phi^a_{-, N}(v_-)``=\Phi^a_{+, N}(v_-)''=J(v_+)\partial_s v_-$ with $J(v_+)$ being considered  as almost complex structure along $v_-.$
To make  sense  out of this senseless identity,   we note that  the  sub-cylinder  with $t\in (-d-l, d+l)$ is  where  the splicing matrix $T_{\beta}$ is not a constant. It is  also  the  place where the loss of differentiability takes place  for the choice of  $\Phi^a_{-, N}$ similar to the one in \cite{2}. Since the sub-cylinder  is contained in both $S^a_-$ and $S^a_+$, we may  define the almost complex structure $J(v_-)$ along $(-d-l, d+l)\times S^1$ to be $J(v_+).$ More generally, we need to define a new extended gluing ${\hat v}_+:{\hat S}_+^a:\rightarrow E$ with ${\hat S}_+^a\simeq (-\infty, \infty)\times S^1\simeq {S}_-^a$ (see the definition in next section). Denote the identification map by $\Gamma:{ S}_-^a\rightarrow {\hat S}_+^a$, which will transfer the almost complex structure $J({\hat v}_+)$ into $J({\hat v}_+\circ \Gamma)$ considered as an almost complex structure along $v_-$.
Let $J(v) =J({\hat v}_+\circ \Gamma)\oplus J(v_+)$  be the corresponding total almost complex structure  along  $v=(v_-, v_+).$  Then $\Phi^a_{ N}$ is defined to be $ \Phi^a_{ N}(v)=J(v)\partial_s v.$   The idea of this construction is to enforce the commutativity of the splicing  matrix $T_{\beta}$ and the total almost complex structure $J(v)$ along $v$: $J(v)\circ T_{\beta}=T_{\beta}\circ J(v).$  One can see  from the proof of the main theorem below that this commutativity is essentially equivalent to the requirement  of no loss of differentiability for $ \Phi_{ N}$.

Then  the main theorem of this paper is the following theorem.

\begin{theorem}
	Using the gluing profile $R=e^{1/r}-e^{1/r_0}$, the filled-section $\Psi_N=\{\Psi^{R_{\theta}}_N\}$ above  with $r\in [0, r_0)$ (hence $R\in (R_0, \infty]$  for $R_{\theta}=(R,\theta)$ and $\theta\in S^1$) is of class $C^1$. Consequently, the filled-section $\Psi=\Psi_L+\Psi_N: L_{k, \delta}^p(C_-, E)\times  L_{k, \delta}^p(C_+, E)\times D_{r_0}\rightarrow L_{k-1, \delta}^p(C_-, E)\times  L_{k-1, \delta}^p(C_+, E)$ is of class $C^1.$
\end{theorem}



Clearly the argument in these sequence of papers can be generalized to deal with  general  nonlinear equations with quasi-linear principal "symbols" as well as higher order equations of similar types. 
Applications of this kind  will be given somewhere else.

\begin{remark}

	(A)	For any given positive integer $m$, the  $C^m$-smoothness  of the filled-section  $\Psi$ can be proved by essentially the same argument but using  different length and center functions $l=L_m(R)=R^{m/(m+1)}\cdot ln^2 R$ and $d=3l=3L_m(R)$. For $C^1$-smoothness here, the function $l=L_1 (R)$. 
	
\medskip
\noindent
	(B) One may assume that $ \Phi_{-, N}$ has the general form  $ \Phi_{ -, N}(v)=A(v_-)\cdot \partial _s v_-$ with $A(v_-) $ being a $End (E)$-valued section over $v_- $ which is to be chosen. Let $J_A(v)=A(v_-)\oplus J(v_+).$
	Then one can show that for generic $J$,  the resulting  $ \Psi_{ N}$ has no loss of differentiability   implies  the  commutativity:$J_A(v)\circ T_{\beta}=T_{\beta}\circ J_A(v)$ so that 
	upto the choices of the transfer map $\Gamma$, our definition for $ \Psi_{-, N}$ is essentially the only possible choice with the desired properties.
	
\end{remark}
The main theorem will be proved in Sec. 3 after recalling the basic definitions in the splicing in Sec. 2.

Like {\cite 3}, we will only deal with the case for $L_{k, \delta}^p$-maps with fixed ends at their double points. Throughout this paper we will assume that (i) $p>2$ and (ii) $k-2/p>1$. Unlike \cite {3}, this last condition is needed for the
 estimates in Sec. 3 in order to prove the main theorem above. If one insists using Hilbert space  (hence $p=2$), then the condition above becomes $k\geq 3.$

\section{Basic definitions of the  splicing}
In this section we recall the definitions of the splicing in \cite{3}.
These definitions are tailored for defining the filled-section in the setting of Banach analysis (comparing with the definitions in \cite{2}).

\subsection{Total gluing of the nodal surface  $S$}


Let $S=C_-\cup_{d_-=d_+}C_+$   with the  double point
${d_-=d_+}$ with  $(C_{\pm}, d_{\pm})$ being  the standard disk  $ (D, 0)$.
Identify $(C_{-}, d_{-})$ with $ ((-\infty, 0)\times S^1, -\infty\times S^1 )=(L_-\times S^1, -\infty\times S^1)$
canonically upto a rotation by  considering the double point $d_-$ as the $S^1$ at $-\infty$ of the half cylinder $L_-\times S^1$.   Here we have denoted the negative half line $(-\infty, 0)$ by $L_-$.
Similarly $(C_{+}, d_{+})\simeq (( 0, \infty,)\times S^1, \infty\times S^1 )=(L_+\times S^1, \infty\times S^1).$

\medskip
\noindent ${\bf \bullet}$  {\bf Cylindrical coordinates on $C_{\pm}:$}

\medskip
\noindent
By the identification  $C_{\pm}\simeq L_{\pm}\times S^1$,   each $C_{\pm}$ has  the  cylindrical coordinates $(t_{\pm}, s_{\pm})\in L_{\pm}\times S^1.$

Let $a=(R, \theta)\in [0, \infty]\times S^1$ be the gluing parameter. 
To defined the total gluing/deformation $S^a=S^{(R, \theta)}$ with gluing parameter $R\not= \infty$, we  introduce  the $a$-dependent cylindrical coordinates $(t^{a, {\pm}}, s^{a,{\pm}})$ on  $C_{\pm}$ by the formula 
$t_{\pm}=t^{a, \pm }\pm  R$ and $s_{\pm }=s^{a, \pm}\pm\theta$. In the following if there is no confusion, we will  denote  $t^{a, \pm }$ by $t$ and $s^{a, \pm}$ by $s$ for both of these $a$ -dependent cylindrical coordinates.

Thus the $t$-range for $L_-$ is $(-\infty, R)$ and  the $t$-range for $L_+$ is $(-R, \infty)$ with the intersection $L_-\cap L_+=(-R, R)$.

\medskip
\noindent ${\bf \bullet}$  {\bf Total gluing $S^a=(S^a_-, S^a_+)$ :}

\medskip
\noindent

In term of the $a$-dependent cylindrical coordinates $(t,s)$,   $C_-=(-\infty, R)\times S^1$ and $C_+=(-R,\infty)\times S^1$.

Then  $S^a_+$ is defined to be the finite cylinder of length $2R$  obtained by gluing $(-1, R)\times S^1\subset C_-$
with $(-R, 1)\times S^1\subset C_+$ along the "common" region $(-1, 1)\times S^1$ by the identity map in term of he $a$-dependent  coordinates $(t,s)$.
Similarly,  $S^a_-$ is the infinite cylinder defined by  gluing $(-\infty, 1)\times S^1\subset C_-$
with $(-1, \infty)\times S^1\subset C_+$ along  $(-1, 1)\times S^1$ by the identity map.

Geometrically, both $S^a_{\pm}$ are obtained by first cutting  each $C_{\pm}$ along the circle at $t=0$ into two sub-cylinders,then gluing back  the  sub-cylinders in $C_-$ with  the corresponding ones in $C_+$ along  the same circle  with a relative rotation of angle $2\theta$.
Set $S^{\infty}=S.$

Now the  cylindrical coordinates $(t_{\pm}, s_{\pm})$ on $C_{\pm}$ as well as the $a$-dependent cylindrical coordinates
$(t, s)$  become the corresponding ones on each  $S^{a}_\pm$ with the relation:   $t_\pm=t\pm R$ and $s_\pm =s\pm \theta.$






%
\subsection{ Splicing matrix $T_{\beta}$}
To defined   $T_{\beta}$,   we need to choose a
length function depending $R$. For the purpose of this paper, the length fuction  $L(R)=L_1(R)=R^{1/2} \cdot ln^2 R$. 

The  splicing matrix $T_{\beta}$ used in this paper is defined by using a  pair of  cut-off function $\beta=(\beta_-,\beta_+)$ depending on the two  parameters $(l, d)$ that parametrize the group of affine transformations $\{t\rightarrow lt+d\}$ of ${\bf R}^1$ defined as follows.

Fix a smooth  cut-off function $\alpha_-:{\bf R}^1\rightarrow [0, 1]$ with the property
that  $\alpha_-  (t)=1$ for $ t<-1$,  $\alpha (t)=0 $ if $ t>1 $ and $\alpha'\leq 0$.
Let  $\alpha_+=:1-\alpha$.
Fix  $l_0>1$ and $d_0>1$. Then  $\beta_{\pm}=\{\beta_{\pm;l,d}\}:{\bf R}^1\times [l_0, \infty)\times [d_0, \infty)\rightarrow [0, 1]$   defined by $\beta_{\pm }(t, l,d)=\alpha_{\pm}( (t\pm d)/l)$, or $\beta_{\pm ;l,d}= \rho_l\circ \tau_{\pm d} \alpha$. Here the translation and multiplication operators are  defined   by $\tau _{d} (\xi)(t)=\xi(t+d)$ and $\rho_l(\xi)(t)=\xi(t/l)$ respectively.

The pair $(l, d)$ will  be  the functions on $R$, $(l=L(R),d=d(R)) $  with   $d= 3l$ and $l$ defined above.

Clearly $\beta_{\pm}$ is a smooth cut-off function  with the following two properties:

\medskip
\noindent
$P_1:$ for  $k\leq k_0$ the $C^0$-norm of the $k$-th derivative $\|\beta_{\pm}^{(k)}\|_{C^0}\leq C/l^k$  , where $C=\|\alpha\|_{C^{k_0}};$

\medskip
\noindent
$P_2:$ under the assumption that $d\geq 3l$, the support of $\beta_{-}'$  is contained in the interval$(d-l,  d+l)$ with  $\beta_{-}=1 $ on   $(-\infty,d-l]$ and $\beta_{-}=0 $ on   $[d+l, \infty);$ and $\beta_{+}'$  is contained in the interval$(-d-l,  -d+l)$ with  $\beta_{+}=1 $ on   $[-d+l, \infty)$ and $\beta_{+}=0 $ on   $(-\infty, -d-l]$.

The splicing matrix then is  defined by

$$T_{\beta}=\left[\begin{array}{ll}
\beta_- &  -\beta_+\\
\beta_+   & \beta_-
\end{array}\right]. $$

Note that from $P_2$, on $(-d+l, d-l)$, $\beta_-=\beta_+=1$.
Then  for $t$  in the  three intervals
$(-\infty, -d-l)$,$(-d+l, d-l)$ and $(d+l, \infty)$,  $T_{\beta}(t)$    are the following constant matrices

$$M_1=Id=\left[\begin{array}{ll}
1 &  0\\
0   &1
\end{array}\right], M_2=\left[\begin{array}{ll}
1&  -1\\
1  & 1
\end{array}\right], \, \, and\,\, M_3= \left[\begin{array}{ll}
0 & -1\\
1  & 0
\end{array}\right] .$$



Note that  $\beta_{\pm}(t)<1$ implies that $\beta_{\mp}(t)=1$ so that

the determinant   $$1\leq D= det\, \left[\begin{array}{ll}
\beta_-  & -\beta_+ \\
\beta_+ & \beta
\end{array}\right]   =\beta_-^2+\beta_+^2\leq 2.$$

This implies that $T^a= (\ominus_a, \oplus_a)$ defined below is invertible uniformly.

\subsection{Total gluing  $T^a$ of maps and sections}


Let $C^{\infty }(C_{\pm}, E )$ be the set of $E$-valued $C^{\infty}$ functions 
on $C_{\pm}$, where $E={\bf C} ^n.$  Similarly $C^{\infty }(S^a_{\pm}, E)$
consists of all  $E$-valued smooth functions on $S^a_{\pm}$

Then $T^a=(T^a_-, T^a_+):  C^{\infty}(C_{-} E)\times C^{\infty} (C_{+}, E)\rightarrow C^{\infty}(S^a_{-} E)\times C^{\infty}(S^a_{+}, E)$ is defined as follows.

In matrix notation, for each $(\xi_-, \xi_+)\in  C^{\infty}(C_{-}, E)\times C^{\infty}(C_{+}, E)$ considered as a column vector,
$$ T^a((\xi_{-}, \xi_{+}))=(T^a_-(\xi_{-}, \xi_{+}), T^a_+(\xi_{-}, \xi_{+}))=(\xi_{-}\ominus_a\xi_{+},  \xi_{-}\oplus \oplus_a \xi_{+})$$ $$ =\left[\begin{array}{ll}
\beta_- &  -\beta_+\\
\beta_+   & \beta_-
\end{array}\right]
\left[\begin{array}{c}
\tau_{-a}\xi_{-} \\
\tau _{a} \xi_{+}
\end{array}\right]
.$$

The inverse of the total gluing,  $(T^a)^{-1}=(T^a_-, T^a_+)^{-1}:C^{\infty}(S^a_{-}, E)\times C^{\infty}(S^a_{+}, E)  \rightarrow C^{\infty}(C_{-} E)\times C^{\infty}(C_{+}, E)$ is defined as following:
for a pair of the $E$-valued  functions $(\eta_{-}, \eta_{+})\in C^{\infty}(S^a_{-} E)\times C^{\infty}(S^a_{+}, E)$,

$(T^a)^{-1}(\eta_{-}, \eta_{+})=(\oplus_a \oplus \ominus_a)^{-1}(\eta_{-}, \eta_{+})=$
$$\left[\begin{array}{ll}
\tau_{a} & 0 \\
0  &  \tau_{-a}
\end{array}\right] \cdot{\frac {1}{D}}\cdot
\left[\begin{array}{ll}
\beta_- &  \beta_+\\
-\beta_+  & \beta_-
\end{array}\right]
\left[\begin{array}{c}
\eta_{-} \\
\eta _{+}

\end{array}\right].$$




Note that the map $(u_-, u_+)\in C^{\infty}(L_{-}\times S^1,  E)\times C^{\infty}(L_{+}\times S^1, E)$ satisfies the asymptotic condition that
$u_-(-\infty)=u_+(\infty)$ that is corresponding to the condition that $T^a_-(u_-, u_+)(-\infty)=-T^a_-(u_-, u_+)(+\infty)$.
In other words $T^a$ maps these two subspaces each other isomorphically as before.

\subsection{ Definition of the transfer map $\Gamma^R$}

First fix the length function $l=L(R)={\sqrt R} (ln R)^2$ and the center function $d= 3l $ with $R\geq R_0$ for a fixed $R_0>>1.$

${\bullet } $  Extended gluing:  
For the gluing parameter $a=(R, \theta),$
extended gluing $u_-{\hat {\oplus}}_au_+:{\hat S}^a_+:={\bf R}^1\times S^1\rightarrow E$ is defined by :
$$u_-{\hat {\oplus}}_au_+=\gamma_{-, -d}\gamma_{+, d}
u_-{ {\oplus}}_au_+ + \{(1-\gamma_{-, -d})\tau_{-a}u_-+(1-\gamma_{+, d})\tau_{a}u_+,$$ where $\gamma_{+, d}=\tau_{-(d+l)}{\hat \gamma}_+$ and $\gamma_{-, -d}=\tau_{(d+l)}{\hat\gamma}_-$  with ${\hat \gamma}_+(t)={\hat \gamma}_-(-t)$.
Here ${\hat \gamma}_+:{\bf R}^1\rightarrow [0,1 ]$ is a smooth cut-off function such  that
${\hat \gamma}_+(t)=1$ for $t<1$ and  ${\hat \gamma}_+(t)=0$ for $t>2$.

The transfer map $\Gamma^a:S_-^a={\bf R}^1\times S^1\rightarrow  {\hat S}_+^a ={\bf R}^1\times S^1$ is defined to be the identity map in term of the $(t,s)$-coordinates for both $S_-^a$ and ${\hat S}_+^a$.

In the following  $u_-{\hat {\oplus}}_au_+$ will be denoted by
${\hat v}^a_+$ and  $u_-{{\oplus}}_au_+$ by ${ v}^a_+.$

Thus in $t$-coordinate, for $-d-l-1<t<d+l+1$, ${\hat v}^a_+={ v}^a_+$,
for $t>d+l+2$, ${\hat v}^a_+=\tau_a u_+ $ and for $t<-d-l-2$, ${\hat v}^a_+=\tau_{-a}u_-. $


Although, ${\Gamma^a}$ does not has effect on the total gluing $T^a$, it transfers  the complex structure
$J({\hat v}_+)$ along the map $ {\hat v}_+:{\hat S}_+^a\rightarrow M$ to a complex structure $J({\hat v}_+\circ \Gamma^a)$  along  $v_-:S_-^a\rightarrow M$ by composing with the  map $\Gamma^a:S_-^a\rightarrow {\hat S}_+^a$.

Then $\Phi^a_+(v_+)={\partial_tv_+}+J(v_+)({\partial_s v_+})$ is  defined to be  the usual
${\overline{\partial}_J}$-operator.
The new $\Phi^a_-$ is  still  a
${\overline{\partial}}$-operator  but with the above almost complex structure
on $v_-$  : $\Phi^a_-(v_-)={\partial_t v_-}+\omega_-\cdot (v_-, v_+)+J({\hat v}^a_+\circ {\Gamma^a})({\partial_s v_-})$. Here $\omega=(\omega_-, \omega_+)$ is the connection matrix introduced in \cite {3}. 

We note that in order to makes sense for the definition of $\Phi^a_-$, we need to know that ${\hat v}^a_+$ is determined by ${ v}^a_{\pm}.$ Indeed, given ${ v}^a_{\pm}$, there is an unique pair  $(u_-, u_+) (=(T^a)^{-1}(v^a_-, v^a_+))$ such that  $T^a(u_-, u_+)=(v^a_-, v^a_+)$. Then we can construct the extended gluing ${\hat v}^a_+=:u_-{\hat {\oplus}}_au_+$ as above.

Let $\Phi^a=(\Phi^a_-, \Phi^a_+)$.
Then $\Psi^a=:(T^a)^{-1}\circ\Phi^a\circ T^a.$

\section{ The nonlinear part $N$ of  $\Psi$} 

Since the splicing matrix $T_\beta$ is $s$-independent and  the translation operator $\tau_{\theta}$ appeared in the total gluing map $T^a$ commutes with both ${\partial}_t$ and ${\partial }_s$,
$\tau_{\theta}$ does not affect  analysis here in any essential way. In  the most part of the rest of this section we will only give the details for  the results    using $T^R$ and state the corresponding ones using $T^a.$

We now derive  the formula for the nonlinear part  $\Psi^R_N$.

\medskip
\noindent
{\bf The nonlinear part $N^R$:}

The nonlinear part of $\Psi^R$
denoted by $N^R: C^{\infty}(C_-, E)\times C^{\infty}(C_+, E) \rightarrow  C^{\infty}(C_-, E)\times C^{\infty}(C_+, E), $
is defined as follows.

For $ (u_-, u_+)\in C^{\infty}(C_-, E)\times C^{\infty}(C_+, E)$ with $T^R(u_-, u_+)=(v_-, v_+),$
$N^R(u_-, u_+)=(T^R)^{-1}( J(v_+\circ \Gamma^R)\partial_s v_-, J(v_+) \partial_s v_+).$

In matrix notation
$N^R((u_-, u_+)=$
$$\left[\begin{array}{ll}
\tau_{R} & 0 \\
0  &  \tau_{-R}
\end{array}\right] \cdot{\frac {1}{D}}\cdot
\left[\begin{array}{ll}
\beta_- &  \beta_+\\
-\beta_+  & \beta_-
\end{array}\right]
\left[\begin{array}{ll}
J({\hat  v}^R_+\circ \Gamma^R) & 0 \\
0  & J( v^R_+)
\end{array}\right]
$$
$$
\partial_s \left\{ \left[\begin{array}{ll}
\beta_-  &  -\beta_+\\
\beta_+  & \beta_-
\end{array}\right]
\left[\begin{array}{c}
\tau_{-R} {u}_{-} \\
{\tau_{R}	u _{+}}

\end{array}\right] \right\} $$

$$=\left[\begin{array}{ll}
\tau_{R} & 0 \\
0  &  \tau_{-R}
\end{array}\right] \cdot{\frac {1}{D}}\cdot
\left[\begin{array}{ll}
\beta_- &  \beta_+\\
-\beta_+  & \beta_-
\end{array}\right]
\left[\begin{array}{ll}
J( {\hat  v}^R_+\circ \Gamma^R) & 0 \\
0  & J( v^R_+)
\end{array}\right]
$$
$$
\left[\begin{array}{ll}
\beta_- &  -\beta_+\\
\beta_+ & \beta_-
\end{array}\right]
\left[\begin{array}{c}
\partial_s \tau_{-R} {u}_{-} \\
\partial_s {\tau_{R}	u _{+}}
\end{array}\right] . $$


Here $ v^R_+= \beta_+\tau_{-R} {u}_{-}+\beta_-\tau_{R} {u}_{+}.$

We note that starting from the last term of the last identity,  the result of each
matrix multiplication can be interpreted  as a pair of function defined on the common domain ${\bf R}^1$ in $t$-coordinate even they have different domains. This follows from the way that  the total gluing map $T^R$ and its inverse are defined. It can also be seen from our next computation of $N^R((u_-, u_+)$
by restricting the pair of functions to three type of subintervals in $t$-variable.

\medskip
\noindent
$\bullet$ {\bf The $N^R$ over the interval  $(-d-l, d+l)$:}

On the $t$-interval $(-d-l, d+l)$ where the non-trivial part of $\beta$ is lying on ${\hat v }^R_+=v_+^R$ and in $t$-coordinate
$\Gamma^R$ is the identity map so that $J({\hat v }^R_+\circ \Gamma^R)=J(v^+).$ Then all the terms of the following identity are well-defined functions on $(-d-l, d+l)\times S^1$ in $(t,s)$ coordinate, and we have that on $(-d-l, d+l)\times S^1$

$${\frac {1}{D}}\cdot
\left[\begin{array}{ll}
\beta _-&  \beta_+\\
-\beta_+  & \beta_-
\end{array}\right]
\left[\begin{array}{ll}
J( {\hat v }^R_+\circ \Gamma^R) & 0 \\
0  & J( v^R_+)
\end{array}\right]
$$
$$
\left[\begin{array}{ll}
\beta_-  &  -\beta_+\\
\beta_+ & \beta_-
\end{array}\right]
\left[\begin{array}{c}
\partial_s \tau_{-R} {u}_{-} \\
\partial_s {\tau_{R}	u _{+}}
\end{array}\right]  $$

$$={\frac {1}{D}}\cdot
\left[\begin{array}{ll}
\beta_- &  \beta_+\\
-\beta_+  & \beta_-
\end{array}\right]
\left[\begin{array}{ll}
J( v^R_+) & 0 \\
0  & J( v^R_+)
\end{array}\right]
$$
$$
\left[\begin{array}{ll}
\beta_-  &  -\beta_+\\
\beta_+ & \beta_-
\end{array}\right]
\left[\begin{array}{c}
\partial_s \tau_{-R} {u}_{-} \\
\partial_s {\tau_{R}	u _{+}}
\end{array}\right]  .$$
$$={\frac {1}{D}}\cdot
\left[\begin{array}{ll}
\beta_- &  \beta_+\\
-\beta_+  & \beta_-
\end{array}\right]
$$
$$
\left[\begin{array}{ll}
\beta_- J( v^R_+) &  -\beta_+J( v^R_+)\\
\beta_+J( v^R_+) & \beta_- J( v^R_+)
\end{array}\right]
\left[\begin{array}{c}
\partial_s \tau_{-R} {u}_{-} \\
\partial_s {\tau_{R}	u _{+}}

\end{array}\right]  $$

$$={\frac {1}{D}}\cdot
\left[\begin{array}{ll}
\beta_- &  \beta_+\\
-\beta_+  & \beta_-
\end{array}\right]
$$
$$
\left[\begin{array}{ll}
\beta_-  &  -\beta_+\\
\beta_+ & \beta_-
\end{array}\right]
\left[\begin{array}{c}
J( v^R_+) \partial_s \tau_{-R} {u}_{-} \\
J( v^R_+) \partial_s {\tau_{R}	u _{+}}
\end{array}\right]  =\left[\begin{array}{c}
J( v^R_+) \partial_s \tau_{-R} {u}_{-} \\
J( v^R_+) \partial_s {\tau_{R}	u _{+}}
\end{array}\right].$$

Thus
$$N^R((u_-, u_+)=
\left[\begin{array}{ll}
\tau_{R} & 0 \\
0  &  \tau_{-R}
\end{array}\right]\left[\begin{array}{c}
J( v^R_+) \partial_s \tau_{-R} {u}_{-} \\
J( v^R_+) \partial_s {\tau_{R}	u _{+}}
\end{array}\right]=\left[\begin{array}{c}
\tau_{R} J( v^R_+) \partial_s {u}_{-} \\
{\tau_{-R}}J( v^R_+) \partial_s 	u _{+}
\end{array}\right].$$

$N^R((u_-, u_+)=(\tau_{R} J( v^R_+) \partial_s {u}_{-},  {\tau_{-R}}J( v^R_+) \partial_s 	u _{+})$

\medskip
\noindent
$\bullet$ {\bf The $N^R$  away from the interval $(-d-l+1,d+ l-1)$}:

We already know from last subsection that away from $(-d-l,d+ l)$,
$T_{\beta}$ is given by $M_1=id$ and $M_3$. We may assume that this is true away from
the interval $(-d-l+1,d+ l-1). $
Hence
for $t>d+l-1$,
$D=1$

and $${\frac {1}{D}}\cdot
\left[\begin{array}{ll}
\beta_- &  \beta_+\\
-\beta_+  & \beta_-
\end{array}\right]
\left[\begin{array}{ll}
J( {\hat v }^R_+\circ \Gamma^R) & 0 \\
0  & J( v^R_+)
\end{array}\right]
$$
$$
\left[\begin{array}{ll}
\beta_-  &  -\beta_+\\
\beta_+ & \beta_-
\end{array}\right]
\left[\begin{array}{c}
\partial_s \tau_{-R} {u}_{-} \\
\partial_s {\tau_{R}	u _{+}}
\end{array}\right]  $$
$$=
\left[\begin{array}{ll}
0 &  1\\
-1  & 0
\end{array}\right]
\left[\begin{array}{ll}
J( {\hat v }^R_+\circ \Gamma^R) & 0 \\
0  & J( v^R_+)
\end{array}\right]
\left[\begin{array}{ll}
0  &  -1\\
1  & 0
\end{array}\right]
\left[\begin{array}{c}
\partial_s \tau_{-R} {u}_{-} \\
\partial_s {\tau_{R}	u _{+}}
\end{array}\right]  $$

$$=
\left[\begin{array}{ll}
0 &  1\\
- 1  & 0
\end{array}\right]
\left[\begin{array}{ll}
J( {\hat v }^R_+\circ \Gamma^R) & 0 \\
0  & J( v^R_+)
\end{array}\right]
\left[\begin{array}{c}
- \partial_s \tau_{R} {u}_{+} \\
\partial_s {\tau_{-R}	u _{-}}
\end{array}\right]  $$

$$=
\left[\begin{array}{ll}
0 &  1\\
- 1  & 0
\end{array}\right]
\left[\begin{array}{c}
- J( {\hat v }^R_+\circ \Gamma^R) \partial_s \tau_{R} {u}_{+} \\
J( v^R_+) \partial_s {\tau_{-R}	u _{-}}
\end{array}\right] =  \left[\begin{array}{c}
J( v^R_+) \partial_s {\tau_{-R}	u _{-}}  \\
J( {\hat v }^R_+\circ \Gamma^R) \partial_s \tau_{R} {u}_{+}
\end{array}\right] .$$

Note that in above computation, each steps makes sense even  the pair of functions are not defined  on the same domain.

Hence for $t>d+l-1$,
$N^R (u_-, u_+) =(\tau_{R} J( v^R_+) \partial_s {	u _{-}}, \tau_{-R}  J( {\hat v}^R_+\circ \Gamma^R) \partial_s {u}_{+}).$

For $t<-d-l+1$, both matrices become the identity matrix with $D=1$ so that

$${\frac {1}{D}}\cdot
\left[\begin{array}{ll}
\beta_- &\beta_+\\
-\beta_+  & \beta_-
\end{array}\right]
\left[\begin{array}{ll}
J( {\hat v }^R_+\circ \Gamma^R) & 0 \\
0  & J( v^R_+)
\end{array}\right]
$$
$$
\left[\begin{array}{ll}
\beta_-  &  -\beta_+\\
\beta_+ & \beta_-
\end{array}\right]
\left[\begin{array}{c}
\partial_s \tau_{-R} {u}_{-} \\
\partial_s {\tau_{R}	u _{+}}
\end{array}\right]  $$
$$=
\left[\begin{array}{ll}
J({\hat v }^R_+\circ \Gamma^R) & 0 \\
0  & J( v^R_+)
\end{array}\right]
\left[\begin{array}{c}
\partial_s \tau_{-R} {u}_{-} \\
\partial_s {\tau_{R}	u _{+}}
\end{array}\right] = \left[\begin{array}{c}
J( {\hat v }^R_+\circ \Gamma^R) \partial_s \tau_{-R} {u}_{-} \\
J( v^R_+) \partial_s {\tau_{R}	u _{+}}
\end{array}\right].$$

Hence
for $t<-d-l+1,$
$N^R (u_-, u_+) =(\tau_{R} J( {\hat v }^R_+\circ \Gamma^R) \partial_s {	u _{-}}, \tau_{-R}  J( v^R_+) \partial_s {u}_{+}).$

In summary, we have proved

\begin{lemma}
	In $(t, s)$-coordinate, $N^R (u_-, u_+)=(\tau_{R} J( {\hat v }^R_+\circ \Gamma^R) \partial_s {	u _{-}}, \tau_{-R}  J( v^R_+) \partial_s {u}_{+})$ for $t<-d-l+1$,
	$N^R (u_-, u_+)=(\tau_{R} J( v^R_+) \partial_s {u}_{-},  {\tau_{-R}}J( v^R_+) \partial_s 	u _{+})$ for $-d-l<t<d+l$ and $N^R (u_-, u_+)=(\tau_{R} J( v^R_+) \partial_s {	u _{-}}, \tau_{-R}  J( {\hat v }^R_+\circ \Gamma^R) \partial_s {u}_{+})$ for $t>d+l-1.$

\end{lemma}

Since on the overlap regions, $\Gamma^R$ is the identity map, we have the following

\begin{cor}
	On the overlap regions, $N^R$ is defined by the same formula. There is no loss of differentiability in $N^R$.
\end{cor}

We note that in the above lemma the two components  $N^R_{\pm} (u_-, u_+)$ already expressed   their  own natural coordinates $t_{\pm}$ coming from $u_{\pm},$ but  we still use the $t$-coordinate to  divide
the domains.
It is more convenient to the division of the domains for
$N^R (u_-, u_+)=(N^R_-(u_-, u_+), N^R_+(u_-, u_+))$ in the  natural coordinates $t_{\pm}$. By abusing the notations, in next  lemma coordinates $t_{\pm}$  are still denoted by $t$.

Note that for $d +l-1<t<R$, in $t_-$-coordinate with $t_-=t-R$, it is in  $-R+d+l-1<t_-<0,$ which corresponds to the right half the region the usual gluing $u_-\oplus_R u_+$ coming from $u_-$. For such a  $t$,  $ v^R_+= \beta_+\tau_{-R} {u}_{-}+\beta_- \tau_{R} {u}_{+}= \tau_{-R} {u}_{-}$ so that the term $\tau_R J( v^R_+)=J(\tau_R\circ \tau_{-R} {u}_{-})=J(u_-)$. Similarly for $t<-d-l-2$, it is in  $t_-<-R-d-l-2,$
the term ${\hat v}^R_+\circ \Gamma^R =\tau_{-R }u_-$ again so that
$ \tau_{R} J( {\hat v}^R_+\circ \Gamma^R)=(\tau_R\circ \tau_{-R} {u}_{-})=J(u_-)$.

Then  we have 

\begin{lemma}
	
	$$N^R_-(u_-, u_+)=\left \{ \begin{array}{ll}
	\tau_{R} J( {\hat v}^R_+\circ \Gamma^R) \partial_s {	u _{-}},&  if\, t<-R-d-l-2 ,\\\tau_{R} J( v^R_+) \partial_s {u}_{-}  &  if\, -R-d-l-3<t<-R+d+l+1 \\
	\tau_{R} J( v^R_+) \partial_s {u}_{-}  &  if\, -R+d+l<t<0 \\
	\end{array}\right.$$

	$$=
	\left \{ \begin{array}{ll}
	\tau_{R} J( {\hat v}^R_+\circ \Gamma^R) \partial_s {	u _{-}},&  if\, t<-R-d-l-2 ,\\\tau_{R} J( v^R_+) \partial_s {u}_{-}  &  if\, -R-d-l-3<t<-R+d+l+1 \\
	J( u_-) \partial_s {u}_{-}  &  if\, -R+d+l<t<0 \\
	\end{array}\right.$$
	
	$$=
	\left \{ \begin{array}{ll}
	\ J( u_-) \partial_s {	u _{-}},&  if\, t<-R-d-l-2 ,
	\\\tau_{R} J( {\hat v}^R_+\circ \Gamma^R) \partial_s {	u _{-}},&  if\, -R-d-l-3<t<-R-d-l+1 ,\\\tau_{R} J( v^R_+) \partial_s {u}_{-}  &  if\, -R-d-l<t<-R+d+l+1 \\
	J( u_-) \partial_s {u}_{-}  &  if\, -R+d+l<t<0 \\
	\end{array}\right.$$

	Similarly,
	$$N^R_+(u_-, u_+)=
	\left \{ \begin{array}{ll}
	\tau_{-R} J( {\hat v}^R_+\circ \Gamma^R) \partial_s {	u _{+}},&  if\, t>R+d+l+2 ,\\\tau_{-R} J( v^R_+) \partial_s {u}_{+}  &  if\, R-d-l-1<t<R+d+l+3 \\
	J( u_+) \partial_s {u}_{+}  &  if\, 0<t< R-d-l\\
	\end{array}\right.$$
	
	$$=
	\left \{ \begin{array}{ll}
	\ J( u_+) \partial_s {	u _{+}},&  if\, t>R+d+l+2 ,
	\\\tau_{-R} J( {\hat v}^R_+\circ \Gamma^R) \partial_s {	u _{+}},&  if\, R+d+l-1<t<R+d+l+3 ,\\\tau_{-R} J( v^R_+) \partial_s {u}_{+}  &  if\, R-d-l-1<t<R+d+l \\
	J( u_+) \partial_s {u}_{+}  &  if\,0 <t<R-d-l\\
	\end{array}\right.$$

	$$=
	\left \{ \begin{array}{ll}
	\ J( u_+) \partial_s {	u _{+}},&  if\, t>R+d+l+2 ,
	\\\tau_{-R} J( {\hat v}^R_+\circ \Gamma^R) \partial_s {	u _{+}},&  if\, R-d-l-3<t<R+d+l+3 ,\\
	J( u_+) \partial_s {u}_{+}  &  if\,0 <t<R-d-l-2\\
	\end{array}\right.$$
	
\end{lemma}

In the last identity, we rearrange the division of the domain so that away from $ [R-d-l-3, R+d+l+3 ]$ of length $2(d+l+3)$, $N^R_+(u_-, u_+)$ is just $N^{\infty}_+(u_-, u_+)= J( u_+) \partial_s {u}_{+} $
so that the error term $E^R(u_-, u_+)=:N^R_+(u_-, u_+)-N^{\infty}_+(u_-, u_+)=\{\tau_{-
	R} J( {\hat v}^R_+\circ \Gamma^R)-J(u_+)\} \partial_s {	u _{+}}=\{J( {\hat v}^R_+\circ \Gamma^R\circ \tau_{-
	R} )-J(u_+)\} \partial_s {	u _{+}}$ is localized on $[R-d-l-3, R+d+l+3 ]\times S^1.$

The formulas for $N^{R_\theta}_+(u_-, u_+)$ can be obtained similarly by simply replacing $R$ by $R_\theta$ in above formulas.
\begin{theorem}
	$N=\{N^R\}:L_{k,\delta}^p(C_-, E)\times L_{k,\delta}^p(C_-, E)\times [R_0, \infty]\rightarrow L_{k-1,\delta}^p(C_-, E)\times L_{k-1,\delta}^p(C_-, E)$  is  continuous.
\end{theorem}

\proof

From above formula, the continuity of $N$  for $R\not=\infty$ is clear.
To see 
the continuity of $N$   at  $R=\infty$,we only need to consider $N_+$. Then it follows from the estimate proved below in this section that 
for   $t\in [R-d-l-3, R+d+l+3], $ the term
$$\|J({\hat v}^R_+\circ \Gamma^R\circ \tau_{-R})-J(u_+)||_{C^{k-1}}$$
$$ \leq C(\beta, \gamma)\cdot\|J\|_{C^{k-1}}\cdot e^{-\delta R/2}\cdot  ||u\|_{{k, p,\delta}}(1+e^{-\delta R/2}\cdot||u\|_{{k, p,\delta}}).$$

so that
$$\|N^R_+(u_-, u_+)-N^{\infty}_+(u_-, u_+)\|_{{k-1, p,\delta}}=\|\{J( {\hat v}^R_+\circ \Gamma^R\circ \tau_{-
	R} )-J(u_+)\} \partial_s {	u _{+}}\|_{{k-1, p,\delta}}$$ $$\leq \|\{J( {\hat v}^R_+\circ \Gamma^R\circ \tau_{-
	R} )-J(u_+)\}\|_{C^{k-1}} ||
\partial_s {	u _{+}}\|_{{k-1, p,\delta}}$$ $$ \leq C\cdot  e^{-\delta R/2}\cdot  ||u\|^2_{{k, p,\delta}}(1+e^{-\delta R/2}\cdot||u\|_{{k, p,\delta}}).$$  This together with the fact that $N^{\infty}$ is continuous implies that $N$ is continuous. Indeed $$\|N^R_+(u_1)-N^{\infty}_+(u_2)\|_{{k-1, p,\delta}}\leq \|N^R_+(u_1)-N^{\infty}_+(u_1)\|_{{k-1, p,\delta}}+\|N^{\infty}_+(u_1)-N^{\infty}_+(u_2)\|_{{k-1, p,\delta}}$$

$$\leq C\cdot  e^{-\delta R/2}\cdot  ||u_1\|^2_{{k, p,\delta}}(1+e^{-\delta R/2}\cdot||u_1\|_{{k, p,\delta}})+\|N^{\infty}_+(u_1)-N^{\infty}_+(u_2)\|_{{k-1, p,\delta}}, $$

which is less than $$\epsilon +\|N^{\infty}_+(u_1)-N^{\infty}_+(u_2)\|_{{k-1, p,\delta}} $$ with any give $\epsilon >0$  and  any fixed $u_1$ by taking $R$ sufficiently large. The conclusion then follows from the continuity of $N^{\infty}.$
\QED

\begin{theorem}
	$N=\{N^R\}:L_{k,\delta}^p(C_-, E)\times L_{k,\delta}^p(C_-, E)\times (R_0, \infty)\rightarrow L_{k-1,\delta}^p(C_-, E)\times L_{k-1,\delta}^p(C_-, E)$  is of class $C^1$.
\end{theorem}

\begin{theorem}
	The derivative  $DN=\{DN^R\}:L_{k,\delta}^p(C_-, E)\times L_{k,\delta}^p(C_-, E)\times (R_0, \infty)\rightarrow L(L_{k,\delta}^p(C_-, E)\times L_{k,\delta}^p(C_-, E)\times {\bf R}^1, L_{k-1,\delta}^p(C_-, E)\times L_{k-1,\delta}^p(C_-, E)) $  can be  extended continuously over $r=0$ ( hence $R=\infty$).
\end{theorem}

\begin{cor}
	The derivative  $DN=\{DN^{r_{\theta}}\}:L_{k,\delta}^p(C_-, E)\times L_{k,\delta}^p(C_-, E)\times D_{r_0}\rightarrow L(L_{k,\delta}^p(C_-, E)\times L_{k,\delta}^p(C_-, E)\times {\bf R}^2, L_{k-1,\delta}^p(C_-, E)\times L_{k-1,\delta}^p(C_-, E)) $   is of class $C^1.$
\end{cor}

To prove the theorems, we use the  formula  above for
$N^R (u_-, u_+)=(N^R _-(u_-, u_+), N^R_+ (u_-, u_+)_+).$
Since the two parts $N^R$ are in the same natural, we only need to deal with $N^R_+ (u_-, u_+).$

Let $N=\{N^R, R\in [0, \infty]\}:L_{k,\delta}^p(C_-, E)\times L_{k,\delta}^p(C_+, E)\times [R_0, \infty]\rightarrow  L_{k-1,\delta}^p(C_-, E)\times L_{k,\delta}^p(C_+, E).$

Denote $L_{k,\delta}^p(C_\pm, E)$ by $W_{\pm}$ and $W={W}_-\times { W}_+. $

In these notation,  $N_+:{W}\times [0, \infty]\rightarrow   L_{k-1,\delta}^p(C_+, E)$.

${\bf\bullet }$ Partial derivative $D_{W} N_+$:

Consider the case $R\not = \infty$ first.

By the last formula above, $D_{W} N_+:{W}\times [0, \infty)
\rightarrow L(L_{k,\delta}^p(C_-, E)\times L_{k,\delta}^p(C_-, E)\times {\bf R}^1, L_{k-1,\delta}^p(C_+, E))$  along $(u_-, u_+)$-direction  is given by the following formula:

$$(D_{W} N_+)_{u_-, u_+, R}(\xi_-, \xi_+)=$$
$$\left \{ \begin{array}{ll}
\{\partial J_{ {\hat v}^R_+\circ \Gamma^R\circ\tau_{-R}}D_W ({\hat v}^R_+\circ \Gamma^R\circ\tau_{-R})(\xi)\}\partial_s {u}_{+} +\tau_{-R} J( {\hat v}^R_+\circ \Gamma^R) \partial_s {\xi}_{+} \\if\, R-d-l-3<t<R+d+l+3,\\
\partial J_{ u_+} (\xi_+) \partial_s {u}_{+}+J( u_+) \partial_s {\xi}_{+}\\   if\, t\not \in [R-d-l-3,R+d+l+3]
\end{array}\right.$$

Thus away from $[R-d-l-3, R+d+l+3 ]\times S^1,$
$(D_{W} N_+)_{u_-, u_+, R}(\xi_-, \xi_+)= (D_{W} N^{\infty}_+)_{u_-, u_+}(\xi_-, \xi_+)=\partial J_{ u_+} (\xi_+) \partial_s {u}_{+}+J( u_+) \partial_s {\xi}_{+};$

Let $\rho:[0, \infty)\rightarrow [0, 1]$ be a smooth cut-off function
such that $$\rho (t)=
\left \{ \begin{array}{ll}
1 & \, \, if\,\, t \in [R-d-l-3,R+d+l+3],\\
0  &\, \, if\, t\not \in [R-d-l-2,R+d+l+2]
\end{array}\right.$$

Then $D_W N_+^R=D_W N_+^{\infty}+\rho D_WE_+^R$  with $E^R_+=: N_+^R-N^{\infty}_+,$ so that $\rho D_WE^R_+$  becomes a map $\rho D_WE_+: {W}\times [0, \infty)
\rightarrow L(L_{k,\delta}^p(C_-, E)\times L_{k,\delta}^p(C_-, E)\times {\bf R}^1, L_{k-1,\delta}^p(C_+, E)).$

Since $D_W N_+^{\infty}$ is continuous proved, for instance, by Floer,  we only need to consider $\rho D_W E_+^{R}$.

Now we derive a more explicit formula for  $ D_WE.$
To this  end, note that for $  R-d-l-3<t<R+d+l+3$,
$$(D_{W} E_+)_{u_-, u_+, R}(\xi_-, \xi_+) $$ $$=\{\partial J_{ {\hat v}^R_+\circ \Gamma^R\circ \tau_{-
		R}}D_W ({\hat v}^R_+\circ \Gamma^R\circ \tau_{-
	R})(\xi)-\partial J_{u_+}(\xi_+)\}\partial_s {u}_{+} +\{ J( {\hat v}^R_+\circ \Gamma^R\circ \tau_{-
	R}) -J(u_+)\}\partial_s {\xi}_{+}.$$

Recall that ${\hat v}^R_+=u_-{\hat {\oplus}}_Ru_+=\gamma_{-, -d}\gamma_{+, d}
u_-{ {\oplus}}_Ru_+ + (1-\gamma_{-, -d})\tau_{-R}u_-+(1-\gamma_{+, d})\tau_{R}u_+=\gamma_{-, -d}\gamma_{+, d}
(\beta_+\tau_{-R}u_-+\beta_-\tau_{R}u_+) + (1-\gamma_{-, -d})\tau_{-R}u_-+(1-\gamma_{+, d})\tau_{R}u_+ =(\gamma_{- }\gamma_{+}\beta_++(1-\gamma_-))u_-\circ \tau_{-R}+(\gamma_{-}\gamma_{+}\beta_-+(1-\gamma_+))u_+\circ\tau_R.$

In the last identity above we have denoted $\gamma_{\pm, \pm d}$ by
$\gamma_{\pm}.$

Hence  $${\hat v}^R_+\circ \Gamma^R\circ \tau_{-
	R}={\hat v}^R_+\circ \tau_{-
	R}$$$$ =\{(\gamma_{- }\gamma_{+}\beta_++(1-\gamma_-))\circ\tau_{-
	R}\}u_-\circ \tau_{-2R}+\{(\gamma_{-}\gamma_{+}\beta_-+(1-\gamma_+))\circ\tau_{-
	R}\}u_+.$$

\begin{lemma}
	
	$${\hat v}^R_+\circ \Gamma^R\circ \tau_{-
		R}={\hat v}^R_+\circ \tau_{-
		R}$$$$ =\{(\gamma_{- }\gamma_{+}\beta_++(1-\gamma_-))\circ\tau_{-
		R}\}u_-\circ \tau_{-2R}+\{(\gamma_{-}\gamma_{+}\beta_-+(1-\gamma_+))\circ\tau_{-
		R}\}u_+.$$
	
	$$D_W({\hat v}^R_+\circ \Gamma^R\circ \tau_{-
		R})(\xi)$$ $$ =\{(\gamma_{- }\gamma_{+}\beta_++(1-\gamma_-))\circ\tau_{-
		R}\}\xi_-\circ \tau_{-2R}+\{(\gamma_{-}\gamma_{+}\beta_-+(1-\gamma_+))\circ\tau_{-
		R}\}\xi_+.$$
	
	$$D_R({\hat v}^R_+\circ \Gamma^R\circ \tau_{-
		R}) =-2\{(\gamma_{- }\gamma_{+}\beta_+ +(1-\gamma_-))\circ\tau_{-
		R}\}\partial_tu_-\circ \tau_{-2R}$$ $$ +\partial_R\{(\gamma_{- }\gamma_{+}\beta_++(1-\gamma_-))\circ\tau_{-
		R}\}u_-\circ \tau_{-2R}+
	\partial_R \{(\gamma_{-}\gamma_{+}\beta_-+(1-\gamma_+))\circ\tau_{-
		R}\}u_+.$$	
\end{lemma}

\begin{lemma}
	The function $F_1	: {W}\times [0, \infty)
	\rightarrow L(L_{k,\delta}^p(C_-, E)\times L_{k,\delta}^p(C_-, E)\times {\bf R}^1, L_{k-1,\delta}^p(C_+, E))$
	defined by $$F_1(\xi, R)=\rho
	D_W({\hat v}^R_+\circ \Gamma^R\circ \tau_{-
		R})(\xi)$$ $$ =\rho\cdot(\{(\gamma_{- }\gamma_{+}\beta_++(1-\gamma_-))\circ\tau_{-
		R}\}\xi_-\circ \tau_{-2R}+\{(\gamma_{-}\gamma_{+}\beta_-+(1-\gamma_+))\circ\tau_{-
		R}\}\xi_+)$$ is continuous.
	
\end{lemma}

The proof of the lemma follows from a  few facts below that will be used repeatedly:

\medskip
\noindent
(A) $F_+:{W}\times [0, \infty)
\rightarrow L(L_{k,\delta}^p(C_+, E), L_{k-1,\delta}^p(C_+, E))$
defined by
$F(u, R)(\xi)=\tau_{ R}\xi$ is continuous. There is a corresponding function 
$F_-.$

 \medskip
 \noindent
 (B) Any smooth function $f_{\pm }$ such as $f_{\pm }=\beta_{\pm}':C_{\pm }\rightarrow {\bf R}^1$ gives rise
 a $C^{\infty}$-map $F_{ \pm}: {\bf R}^1 \rightarrow  C^m (C_{\pm},{\bf R}^1)$ defined by  $F_{\pm}(R)=f\circ \tau_R$ for any $m$. In particular, we may assume that  $m>>k$. 
 
\medskip
\noindent
(C) Any smooth section such as $J:B\subset M\rightarrow End (E)$ with $B$ being a small ball in $ M$ gives rise
a $C^{\infty}$-map $F_{J, \pm}:  L_{k, \delta}^p (C_{\pm }, B)\rightarrow L_{k, \delta}^p (C_{\pm}, End (E))$ again defined by the composition $u_{\pm}\rightarrow J\circ u_{\pm}.$
 
 \medskip
 \noindent
 (D) The paring
 
 \noindent
 $L_{k}^p(C_{\pm}, E)\times L(  L_{k, \delta }^p(C_{\pm}, E), L_{k, \delta }^p(C_{\pm}, E))\rightarrow L_{k, \delta }^p(C_{\pm}, E)$ is bounded bilinear and hence smooth as long as the space $L_{k }^p(C_{\pm}, E))$ forms  Banach algebra.
 
 \medskip
 \noindent
 (E)
 \noindent
 For $m>>k$, $L(  L_{k, \delta }^p(C_{\pm}, E), L_{k, \delta }^p(C_{\pm}, E))$ is a $C^m(C_{\pm}, {\bf R}^1)$-module and  the multiplication map    $$C^m(C_{\pm}, {\bf R}^1)\times L(  L_{k, \delta }^p(C_{\pm}, E), L_{k, \delta }^p(C_{\pm}, E))\rightarrow L(  L_{k, \delta }^p(C_{\pm}, E), L_{k, \delta }^p(C_{\pm}, E))$$ is bounded bilinear and hence smooth.
 
 The proofs for  (B), (C) and (E) are straightforward and (D)
 is stated  in the first chapter of Lang's book \cite{5} for general Banach spaces.

The property (A) is well-known and was proved for instance in \cite{3}.  In fact, what is needed here is a modified version of (A)  in Lemma 4.4 of  \cite{3}. Given this modified version of (A), the proof of the lemma is almost identical to the corresponding statement in \cite{3}. We leave  the straightforward verification to the readers.

\begin{cor}
	For $R\not = \infty,$ $D_WN$ is continuous.
\end{cor}

To prove that $D_WN$ can be extended continuously over $R=\infty, $
we only need to show this  for
$\rho D_WE_+: {W}\times [0, \infty)
\rightarrow L(L_{k,\delta}^p(C_-, E)\times L_{k,\delta}^p(C_-, E)\times {\bf R}^1, L_{k-1,\delta}^p(C_+, E))$, which follows from
the estimates  below.


 Recall that we only consider the space  $L_{k, \delta}^p(C_{\pm}, E)$ of $L_k^p$-maps that $\delta$-exponentially decay so that the estimates  are  only applicable to the  case with fixed ends. In particular, the values $u_{\pm (d_{\pm})}$ of each map at the double point  is fixed and set to be 
 $0\in E.$

Then we  have
for  $t\in [R-d-l-3, R+d+l+3], $
$${\hat v}^R_+\circ \Gamma^R\circ \tau_{-
	R}(t, s)={\hat v}^R_+\circ \tau_{-
	R}(t, s)$$$$ =\{(\gamma_{- }\gamma_{+}\beta_++(1-\gamma_-))\circ\tau_{-
	R}\}(t, s)u_-(t-2R, s)+\{(\gamma_{-}\gamma_{+}\beta_-+(1-\gamma_+))
\circ\tau_{-
	R}\}(t, s)u_+(t, s)$$ $$=\{(\gamma_{- }\gamma_{+}\beta_++(1-\gamma_-))\circ\tau_{-
	R}\}(t, s)u_-(t_-, s)+\{(\gamma_{-}\gamma_{+}\beta_-+(1-\gamma_+))
\circ\tau_{-
	R}\}(t, s)u_+(t, s) $$  with $t_-\in [-R-d-l-3, -R+d+l+3].$

Similarly
$$D_W({\hat v}^R_+\circ \Gamma^R\circ \tau_{-
	R})(\xi) (t, s)$$ $$ =\{(\gamma_{- }\gamma_{+}\beta_++(1-\gamma_-))\circ\tau_{-
	R}\} (t, s)\cdot \xi_- (t-2R, s)+\{(\gamma_{-}\gamma_{+}\beta_-+(1-\gamma_+))\circ\tau_{-
	R}\} (t, s)\cdot \xi_+ (t, s)$$$$ =\{(\gamma_{- }\gamma_{+}\beta_++(1-\gamma_-))\circ\tau_{-
	R}\} (t, s)\cdot \xi_- (t_-, s)+\{(\gamma_{-}\gamma_{+}\beta_-+(1-\gamma_+))\circ\tau_{-
	R}\} (t, s)\cdot \xi_+ (t, s)$$ with $t_-\in [-R-d-l-3, -R+d+l+3].$

\begin{lemma}
	For  $t\in [R-d-l-3, R+d+l+3], $ 
	
	(I) $$\|\xi_+\|_{C^{k-1}}\leq C\cdot e^{-\delta ([R-d-l-3)}\|\xi_+\|_{k, p,\delta}, $$ with $C=C_k$ independent of $R$.
	
	(II)
	$$  \|{\hat v}^R_+\circ \Gamma^R\circ \tau_{-
		R}\|_{C^{k-1}}\leq C(\beta, \gamma)\cdot e^{-\delta R/2}\cdot  ||u\|_{{k, p,\delta}} .$$	
	
	(III)
	$$\|D_W({\hat v}^R_+\circ \Gamma^R\circ \tau_{-
		R})(\xi)\|_{C^{k-1}}\leq C(\beta, \gamma)\cdot e^{-\delta R/2}\cdot  ||\xi\|_{{k, p,\delta}}.$$ 
	
	(IV)$$\|J({\hat v}^R_+\circ \Gamma^R\circ \tau_{-R})-J(u_+)||_{C^{k-1}}$$
	$$ \leq \leq C(\beta, \gamma)\cdot\|J\|_{C^{k-1}}\cdot e^{-\delta R/2}\cdot  ||u\|_{{k, p,\delta}}(1+e^{-\delta R/2}\cdot||u\|_{{k, p,\delta}}).$$
	
\end{lemma}

\proof

(I) For  $t\in [R-d-l-3, R+d+l+3], $ 
applying Sobolev embedding to each sub-cylinder of length 1 inside $[R-d-l-3, R+d+l+3]$, there is a constant $C=C_k$ independent of $R$ such that 
$$\|\xi_+\|_{C^{k-1}}\leq C\|\xi_+|_{[R-d-l-3, R+d+l+3]}\|_{k, p}\leq e^{-\delta (R-d-l-3)}\|\xi_+\|_{k, p,\delta}.$$

(II) For  $t\in [R-d-l-3, R+d+l+3], $ $$\|{\hat v}^R_+\circ \Gamma^R\circ \tau_{-
	R}\|_{C^{k-1}}\leq C\{1+\|\beta\|_{C^{k-1}}+\|\gamma\|_{C^{k-1}}\}^3\cdot (\|u_-|_{[-R-d-l-3, -R+d+l+3]}||_{C^{k-1}} $$ $$ +||u_+|_{[R-d-l-3, R+d+l+3]}\|_{C^{k-1}})$$ $$
\leq C(||\beta\|_{C^{k-1}}, ||\gamma\|_{C^{k-1}})\cdot e^{-\delta(R-d-l-3)}\cdot ( ||u_-\|_{{k, p,\delta}}+||u_+\|_{{k, p,\delta}})$$ $$\leq C(\beta, \gamma)\cdot e^{-\delta R/2}\cdot ( ||u_-\|_{{k, p,\delta}}+||u_+\|_{{k, p,\delta}})=C(\beta, \gamma)\cdot e^{-\delta R/2}\cdot  ||u\|_{{k, p,\delta}} .$$
Here and below, we  use $ C(\beta, \gamma)$ to denote the  polynomial  function $C(||\beta\|_{C^{k-1}}, ||\gamma\|_{C^{k-1}})$
of positive integer coefficients
in $||\beta\|_{C^{k-1}}$ and  $||\gamma\|_{C^{k-1}}$, which is $C\cdot \{\|\beta\|_{C^{k-1}}+\|\gamma\|_{C^{k-1}}\}^3$ in this case.

(III) Similarly,  for $t\in [R-d-l-3, R+d+l+3], $  the term
$$\|D_W({\hat v}^R_+\circ \Gamma^R\circ \tau_{-
	R})(\xi)\|_{C^{k-1}}\leq C\{1+\|\beta\|_{C^{k-1}}+\|\gamma\|_{C^{k-1}}\}^3\|_{C^{k-1}}\cdot (\|\xi_-|_{[-R-d-l-3, -R+d+l+3]}||_{C^{k-1}} $$ $$ +||\xi_+|_{[R-d-l-3, R+d+l+3]}\|_{C^{k-1}})$$ $$
\leq C(||\beta\|_{C^{k-1}}, ||\gamma\|_{C^{k-1}})\cdot e^{-\delta(R-d-l-3)}\cdot ( ||\xi_-\|_{{k, p,\delta}}+||\xi_+\|_{{k, p,\delta}})$$ $$\leq C(\beta, \gamma)\cdot e^{-\delta R/2}\cdot ( ||\xi_-\|_{{k, p,\delta}}+||\xi_+\|_{{k, p,\delta}}) =C(\beta, \gamma)\cdot e^{-\delta R/2}\cdot  ||\xi\|_{{k, p,\delta}}.$$

(IV) For   $t\in [R-d-l-3, R+d+l+3], $
the term
$$\|J({\hat v}^R_+\circ \Gamma^R\circ \tau_{-R})-J(u_+)||_{C^{k-1}}$$
$$ \leq \|\{J({\hat v}^R_+\circ \Gamma^R\circ \tau_{-R})
-J(u_+)\}|_{[R-d-l-3, R+d+l+3]}||_{C^0}$$ $$ +
\{\|\nabla (J({\hat v}^R_+\circ \Gamma^R\circ\tau_{-R}))|_{[R-d-l-3, R+d+l+3]}||_{C^{k-2}}
+\|\nabla (J((u_+))|_{[R-d-l-3, R+d+l+3]}||_{C^{k-2}}\}$$
$$\leq \|J\|_{C^1}\|({\hat v}^R_+\circ \Gamma^R\circ \tau_{-R}
-u_+)|_{[R-d-l-3, R+d+l+3]}||_{C^0}$$ $$ +
\{\|D J_{{\hat v}^R_+\circ \Gamma^R\circ\tau_{-R}}\circ D({\hat v}^R_+\circ \Gamma^R\circ\tau_{-R}))|_{[R-d-l-3, R+d+l+3]}||_{C^{k-2}}
$$$$ +\|D J_{u_+}\circ D(u_+)|_{[R-d-l-3, R+d+l+3]}||_{C^{k-2}}\}$$

$$\leq \|J\|_{C^1}\{\|{\hat v}^R_+\circ \Gamma^R\circ \tau_{-R}|_{[R-d-l-3, R+d+l+3]}||_{C^0}+||u_+|_{[R-d-l-3, R+d+l+3]}||_{C^0}\}$$ $$ +
\|J\|_{C^{k-1}} \{\| ({\hat v}^R_+\circ \Gamma^R\circ\tau_{-R}))|_{[R-d-l-3, R+d+l+3]}||_{C^{k-1}}^2
$$$$ +\|u_+|_{[R-d-l-3, R+d+l+3]}||_{C^{k-1}}^2\}$$

$$\sim
\|J\|_{C^{k-1}}\cdot $$ $$  \{\| ({\hat v}^R_+\circ \Gamma^R\circ\tau_{-R}))|_{[R-d-l-3, R+d+l+3]}||_{C^{k-1}}(1+\| ({\hat v}^R_+\circ \Gamma^R\circ\tau_{-R}))|_{[R-d-l-3, R+d+l+3]}||_{C^{k-1}})
$$$$ +\|u_+|_{[R-d-l-3, R+d+l+3]}||_{C^{k-1}}(1+\|u_+|_{[R-d-l-3, R+d+l+3]}||_{C^{k-1}})\}$$ $$ \leq C(\beta, \gamma)\cdot\|J\|_{C^{k-1}}\cdot e^{-\delta R/2}\cdot  ||u\|_{{k, p,\delta}}(1+e^{-\delta R/2}\cdot||u\|_{{k, p,\delta}}).$$

\QED

Therefore we have

\begin{pro}
	
	$$\|\rho \cdot (D_{W} E_+)_{u_-, u_+, R}\|_o\sim  e^{-\delta R/2} C(\beta, \gamma)\cdot \|\rho\|_{C^{k-1}} \cdot \| J\|_{C^k}\cdot
	\|u\|_{{k, p, \delta}}\cdot (1 +    ||u\|_{{k, p,\delta}}).
	$$
	
\end{pro}

\proof	

$$\|\rho \cdot (D_{W} E_+)_{u_-, u_+, R}\|_o=:\sup_{\|\xi\|_{k, p, \delta}\leq 1}
\|\rho \cdot (D_{W} E_+)_{u_-, u_+, R}(\xi_-, \xi_+)\|_{k-1, p, \delta} $$
$$=\sup_{\|\xi\|_{k, p, \delta}\leq 1}
\|\rho \cdot\{DJ_{{\hat v}^R_+\circ \Gamma^R\circ \tau_{-  R}}D_W ({\hat v}^R_+\circ \Gamma^R\circ \tau_{-
	R})(\xi)-DJ_{u+}(\xi_+)\}\cdot \partial_s {u}_{+}$$ $$ +\{J({\hat v}^R_+\circ \Gamma^R\circ \tau_{-  R})-J(u_+)\} \cdot\partial_s {\xi}_{+}\|_{k-1, p, \delta}$$

$$ \leq \|\rho\|_{C^{k-1}} \cdot[ \| J\|_{C^k}\cdot
\{\| {\hat v}^R_+\circ \Gamma^R\circ \tau_{-  R}\|_{C^{k-1}}
+ \|u_+\|_{C^{k-1}}\}\cdot $$ $$ \sup_{\|\xi\|_{k, p, \delta}\leq 1}\{\|D_W ({\hat v}^R_+\circ \Gamma^R\circ \tau_{-
	R})(\xi)\|_{C^{k-1}}+\|\xi_+|_{[R-d-l-3, R+d+l+3]}\|_{C^{k-1}}\}\|\partial_s {u}_{+}\|_{k-1, p, \delta}$$ $$  +
\| J({\hat v}^R_+\circ \Gamma^R\circ \tau_{-
	R})-J(u_+)\|_{C^{k-1}} \sup_{\|\xi\|_{k, p, \delta}\leq 1}\|\partial_s {\xi}_{+}\|_{k-1, p, \delta}]$$

(by the estimates in the previous lemma for $\|D_W ({\hat v}^R_+\circ \Gamma^R\circ \tau_{-
	R})(\xi)\|_{C^{k-1}}$ and $\|\xi_+|_{[R-d-l-3, R+d+l+3]}\|_{C^{k-1}}$)
$$ \leq \|\rho\|_{C^{k-1}} \cdot\{ \| J\|_{C^k}\cdot
\{\| {\hat v}^R_+\circ \Gamma^R\circ \tau_{-  R}\|_{C^{k-1}}
+ \|u_+\|_{C^{k-1}}\}$$ $$ \cdot  \sup_{\|\xi\|_{k, p, \delta}\leq 1}\{(1+C(\beta, \gamma))e^{-\delta R/2}\cdot  ||\xi\|_{{k, p,\delta}}\}\|\partial_s {u}_{+}\|_{k-1, p, \delta}
+
\| J({\hat v}^R_+\circ \Gamma^R\circ \tau_{-
	R})-J(u_+)\|_{C^{k-1}} \}$$

(by the estimate in the previous lemma for $\| J({\hat v}^R_+\circ \Gamma^R\circ \tau_{-
	R})-J(u_+)\|_{C^{k-1}} \}$)
$$\sim  \|\rho\|_{C^{k-1}} \cdot\{ \| J\|_{C^k}\cdot
\|u\|_{{k, p, \delta}}^2\cdot (1+C(\beta, \gamma))\cdot   e^{-\delta R/2}
$$ $$ +  C(\beta, \gamma)\cdot\{ \| J\|_{C^k}\cdot e^{-\delta R/2}\cdot  ||u\|_{{k, p,\delta}}(1+||u\|_{{k, p,\delta}})
\}\}$$
$$\sim  e^{-\delta R/2} C(\beta, \gamma)\cdot \|\rho\|_{C^{k-1}} \cdot \| J\|_{C^k}\cdot
\|u\|_{{k, p, \delta}}\cdot (1 +    ||u\|_{{k, p,\delta}}).
$$

\QED

Next consider $D_{R_{\theta}}N_+:L_{k, \delta}^p(C_+, E)\times [R_0, \infty) \rightarrow L_{k-1, \delta}^p(C_+, E)$.

Recall that $$N_+ (R_{\theta},
u_-, u_+) =
\left \{ \begin{array}{ll}
\ J( u_+) \partial_s {	u _{+}},&  if\, t>R+d+l+2 ,
\\\tau_{-R_{\theta}} J( {\hat v}^{R_{\theta}}_+\circ \Gamma^{R_{\theta}}) \partial_s {	u _{+}},&  if\, R-d-l-3<t<R+d+l+3 ,\\
J( u_+) \partial_s {u}_{+}  &  if\,0 <t<R-d-l-2\\
\end{array}\right.$$

Hence  we have

\begin{lemma}
	For $R\not = \infty,$

	$D_{R}N_+ (R,
	u_-, u_+) 
	=$ $$
	\left \{ \begin{array}{ll}
	D J_{ {\hat v}^R_+\circ \Gamma^R\circ\tau_{-R}\circ\tau_{-R}}( D_R\{{\hat v}^R_+\circ \Gamma^R\circ\tau_{-R}\}) \partial_s {	u _{+}},&  if\, t\in [R-d-l-3 ,R+d+l+3 ] ,\\
	0  &  if\, t\not \in [R-d-l-2 ,R+d+l+2]
	\end{array}\right.$$
	
	and  $D_{\theta}N_+ (R_\theta,
	u_-, u_+)  =$ $$
	\left \{ \begin{array}{ll}
	D J_{ {\hat v}^{R_\theta}_+\circ \Gamma^{R_\theta}\circ\tau_{-R_{\theta}}\circ\tau_{-R_{\theta}}}( D_\theta\{{\hat v}^{R_\theta}_+\circ \Gamma^{R_\theta}\circ\tau_{-R_{\theta}}\}) \partial_s {	u _{+}},&  if\, t\in [R-d-l-3 ,R+d+l+3 ] ,\\
	0  &  if\, t\not \in [R-d-l-2 ,R+d+l+2]
	\end{array}\right.$$
\end{lemma}

Recall here

$${\hat v}^R_+\circ \Gamma^R\circ \tau_{-
	R}={\hat v}^R_+\circ \tau_{-
	R}$$$$ =\{(\gamma_{- }\gamma_{+}\beta_++(1-\gamma_-))\circ\tau_{-
	R}\}u_-\circ \tau_{-2R}+\{(\gamma_{-}\gamma_{+}\beta_-+(1-\gamma_+))\circ\tau_{-
	R}\}u_+,$$

and
$$D_R({\hat v}^R_+\circ \Gamma^R\circ \tau_{-
	R}) =-2\{(\gamma_{- }\gamma_{+}\beta_+ +(1-\gamma_-))\circ\tau_{-
	R}\}\partial_tu_-\circ \tau_{-2R}$$ $$ -\{(a_-\gamma_{- }'\gamma_{+}+a_+\gamma_{- }\gamma_{+}'\beta_+-a_-\gamma_-')\circ\tau_{-
	R}\}u_-\circ \tau_{-2R}$$ $$-
\{(a_-\gamma'_{-}\gamma_{+}+a_+\gamma_{-}\gamma_{+}'\beta_--a_+\gamma_+')\circ\tau_{-
	R}\}u_+.$$

Similarly but simpler  $$D_\theta ({\hat v}^{R_\theta}_+\circ \Gamma^{R_\theta}\circ \tau_{-R_\theta
}) =-2\{(\gamma_{- }\gamma_{+}\beta_+ +(1-\gamma_-))\circ\tau_{-
R}\}\partial_su_-\circ \tau_{-2R_\theta}.$$

Here $a_{\pm}=\partial_R (l(R)+d(R))$. Recall that we may assume that
$l(R)+d(R)=R^1/2\cdot \ln^2 R$ for $R\in [R_0, \infty)$ and $R_0>>0$.

Then $|a_{\pm}|\sim R^{-1/2}\cdot \{\ln^2 R-2\ln R\}\sim R^{-1/2}\cdot\ln^2 R.$

From this formula, using the fact that $F: L_{k, \delta}^p(C_+, E)\times [R_0, \infty) \rightarrow L_{k-1, \delta}^p(C_+, E)$ given by $F(u, R)=Du\circ \tau_R$ is continuous together with the general properties labeled as $(B)$, $(C)$ and $(D)$  before,  
it is clear that the next lemma is true.

\begin{lemma}
	For $R\not=\infty,$ $D_RN_{+}$  is continuous
	
\end{lemma}

\begin{pro}
	Assume further that $k-2/p>2$  so that $\|\nabla u\|_{C^0}\leq C \|u\|_{k, p}.$
	Then  $D_RN_{\pm}$ ($D_rN_{\pm}$)extends over  $R=\infty$ ($r=0$) continuously with respect to
	the gluing profile $R=e^{{1/r}}$.  So does  $D_{\theta}N_{\pm}$. Consequently, 
	$N: L_{k, \delta}^p(C_-, E)\times L_{k, \delta}^p(C_+, E)\times D_{r_0}\rightarrow L_{k-1, \delta}^p(C_-, E)\times L_{k-1, \delta}^p(C_+, E)$ is of class $C^1$.
	
\end{pro}

\proof 

This follows from  the following two estimates, in which
the effect of $|a_{\pm}|$ will be ignored since $|a_{\pm}|\sim R^{-1/3}<<1.$

We state each as a lemma.

\begin{lemma}
	For
	$t\in [R-d-l-3 ,R+d+l+3 ]$,
	$$\|D J_{ {\hat v}^R_+\circ \Gamma^R\circ\tau_{-R}\circ\tau_{-R}}( D_R\{{\hat v}^R_+\circ \Gamma^R\circ\tau_{-R}\}) \partial_s {	u _{+}}\|_{k-1, p, \delta}$$
	$$\leq \|DJ\|_{C^{k-1}}\cdot c(\beta, \gamma) \cdot \| {	u}	\|_{k, p, \delta} \cdot(I+II).$$
	Here $I\sim   e^{-\delta R/2} \cdot \|{	u _{+}}	\|_{k, p, \delta} \cdot \| u_- \|_{k, p, \delta}$ and  $II\sim e^{-\delta R/2} \| {	u _{+}}	\|^2_{k, p, \delta} .$
\end{lemma}

\proof

For
$t\in [R-d-l-3 ,R+d+l+3 ]$,

$$\|D J_{ {\hat v}^R_+\circ \Gamma^R\circ\tau_{-R}\circ\tau_{-R}}( D_R\{{\hat v}^R_+\circ \Gamma^R\circ\tau_{-R}\}) \partial_s {	u _{+}}\|_{k-1, p, \delta}$$
$$\leq \|DJ\|_{C^{k-1}}\cdot c(\beta, \gamma) \cdot \| {	u}	\|_{k, p, \delta} \{ \|\partial_s {	u _{+}}	\|_{k-1, p, \delta} \cdot \{\| u_-\circ \tau_{-2R} \|_{C^{0}} +\|\partial_t u_-\circ \tau_{-2R} \|_{C^{0}} \} $$ $$ +\|\partial_s {	 u _{+}}	\|_{k-2, p, \delta} \cdot  \{\| u_-\circ \tau_{-2R} \|_{C^{k-2}} +\|\partial_t u_-\circ \tau_{-2R} \|_{C^{k-2}} \} $$ $$ +\|\partial_s {	u _{+}}	\|_{C^{0}}
\cdot  \{\| u_-\circ \tau_{-2R} \|_{k-1, p, \delta} +\|\partial_t u_-\circ \tau_{-2R} \|_{k-1, p, \delta} \}
$$ $$+ \|DJ\|_{C^{k-1}}\cdot c(\beta, \gamma) \cdot \| {	u}	\|_{k, p, \delta} \cdot \{ \|\partial_s {	u _{+}}	 \|_{k-1, p, \delta} \cdot \| u_+ \|_{C^{0}}   $$ $$ +\|\partial_s {	u _{+}}	\|_{k-2, p, \delta} \cdot  \| u_+ \|_{C^{k-2}}  $$ $$ +\|\partial_s {	u _{+}}	\|_{C^{0}}
\cdot \| u_+\ \|_{k-1, p, \delta} \}=\|DJ\|_{C^{k-1}}\cdot c(\beta, \gamma) \cdot \| {	u}	\|_{k, p, \delta} \cdot(I+II).$$

Estimate for $I$ with $t\in [R-d-l-3 ,R+d+l+3 ]$:

$$I\sim  \|\partial_s {	u _{+}}	\|_{k-1, p, \delta}  \cdot \| u_-\circ \tau_{-2R} |_{[R-d-l-3 ,R+d+l+3 ]}\|_{C^{1}}  $$ $$ +\|\partial_s {	u _{+}}	\|_{k-2, p, \delta} \cdot  \| u_-\circ \tau_{-2R} |_{[R-d-l-3 ,R+d+l+3 ]} \|_{C^{k-1}}   +\|\partial_s {	u _{+}}	 |_{[R-d-l-3 ,R+d+l+3 ]}\|_{C^{0}}
\cdot  \| u_-\circ \tau_{-2R} \|_{k, p, \delta}$$
$$\sim  \|\partial_s {	u _{+}}	\|_{k-1, p, \delta} \cdot  \| u_-\circ \tau_{-2R} |_{[R-d-l-3 ,R+d+l+3 ]}\|_{k, p}  $$ $$ +\|\partial_s {	u _{+}}	\|_{k-2, p, \delta} \cdot  \| u_-\circ \tau_{-2R} |_{[R-d-l-3 ,R+d+l+3 ]} \|_{k, p}
$$ $$ +\|\partial_s {	u _{+}}	 |_{[R-d-l-3 ,R+d+l+3 ]}\|_{C^{0}}
\cdot  \| u_-\circ \tau_{-2R} |_{[R-d-l-3 ,R+d+l+3 ]}\|_{k, p, \delta}
$$

$$\leq  \|{	u _{+}}	\|_{k, p, \delta}  \cdot  \| u_- |_{[-R-d-l-3 ,-R+d+l+3 ]}\|_{k, p}  $$  $$ +\|\partial_s {	u _{+}}	 |_{[R-d-l-3 ,R+d+l+3 ]}\|_{C^{0}}
\cdot  e^{\delta (R+d+l+3)} \| u_-\circ \tau_{-2R} |_{[R-d-l-3 ,R+d+l+3 ]}\|_{k, p}
$$

$$\leq  \|{	u _{+}}	\|_{k, p, \delta}  \cdot e^{-\delta (R-d-l-3)} \|e_-(t) u_- |_{[-R-d-l-3 ,-R+d+l+3 ]}\|_{k, p}  $$  $$ +\|\partial_s {	u _{+}}	 |_{[R-d-l-3 ,R+d+l+3 ]}\|_{C^{0}}
\cdot  e^{\delta (R+d+l+3)} \| u_- |_{[-R-d-l-3 ,-R+d+l+3 ]}\|_{k, p}
$$

$$\leq  \|{	u _{+}}	\|_{k, p, \delta}  \cdot e^{-\delta (R-d-l-3)} \| u_- \|_{k, p, \delta}   +e^{-\delta (R-d-l-3)}\|e_+(t)\partial_s {	u _{+}}	 |_{[R-d-l-3 ,R+d+l+3 ]}\|_{1, p}
$$ $$ \cdot  e^{\delta (R+d+l+3)}\cdot  e^{-\delta (R-d-l-3)}\| e_-(t)u_- |_{[-R-d-l-3 ,-R+d+l+3 ]}\|_{k, p}
$$ $$ \sim   e^{-\delta R/2} \cdot \|{	u _{+}}	\|_{k, p, \delta} \cdot \| u_- \|_{k, p, \delta}   +e^{-\delta (R-d-l-3)}\| {	u _{+}}	 \|_{k, p, \delta}
\cdot  e^{\delta (R+d+l+3)}\cdot  e^{-\delta (R-d-l-3)}\| u_- \|_{k, p}
$$ $$ \sim   e^{-\delta R/2} \cdot \|{	u _{+}}	\|_{k, p, \delta} \cdot \| u_- \|_{k, p, \delta}.$$

Estimate for $II$ with $t\in [R-d-l-3 ,R+d+l+3 ]$:

$$II\sim \| {	u _{+}}	\|_{k, p, \delta}  \cdot  \| u_+ \|_{C^{k-2}} \sim \| {	u _{+}}	\|_{k, p, \delta}  \cdot  \| u_+ |_{[R-d-l-3 ,R+d+l+3 ]}\|_{k-1, p}  $$
$$ \leq \| {	u _{+}}	\|_{k, p, \delta}  \cdot  e^{-\delta(R-d-l-3)}\| e_+(t)u_+ |_{[R-d-l-3 ,R+d+l+3 ]}\|_{k-1, p}   \leq  e^{-\delta R/2} \| {	u _{+}}	\|^2_{k, p, \delta}  .$$

\QED

It follows from the estimate in above lemma that $ D_R N_+$ extends over $R=\infty.$

Finally we need to show that $D_\theta N_+$ extends over $R=\infty$. This follows   from  the following estimate.

\begin{lemma}
	For
	$t\in [R-d-l-3 ,R+d+l+3 ],$ $$ ||D_{\theta}N_+ (R_\theta,
	u_-, u_+) ||_{k-1, p, \delta}\leq\|DJ\|_{C^{k-1}}\cdot c(\beta, \gamma) \cdot \| {	u}	\|_{k, p, \delta}\cdot III$$	with $III\sim e^{-\delta R/2}\cdot \| {	u _{-}}	\|_{k, p, \delta} \cdot\| {	u _{+}}	\|_{k, p, \delta}.$
\end{lemma}

\proof

For
$t\in [R-d-l-3 ,R+d+l+3 ]$,
recall $$D_\theta ({\hat v}^{R_\theta}_+\circ \Gamma^{R_\theta}\circ \tau_{-R_\theta
}) =-2\{(\gamma_{- }\gamma_{+}\beta_+ +(1-\gamma_-))\circ\tau_{-
R}\}\partial_su_-\circ \tau_{-2R_\theta}.$$

Then $||D_{\theta}N_+ (R_\theta,
u_-, u_+) ||_{k-1, p, \delta} =$
$$\|D J_{ {\hat v}^{R_\theta}_+\circ \Gamma^{R_\theta}\circ\tau_{-R_{\theta}}\circ\tau_{-R_{\theta}}}( D_\theta\{{\hat v}^{R_\theta}_+\circ \Gamma^{R_\theta}\circ\tau_{-R_{\theta}}\}) \partial_s {	u _{+}}\|_{k-1, p, \delta}$$
$$\leq \|DJ\|_{C^{k-1}}\cdot c(\beta, \gamma) \cdot \| {	u}	\|_{k, p, \delta} \{ \|\partial_s {	u _{+}}	\|_{k-1, p, \delta} \cdot \|\partial_s u_-\circ \tau_{-2R_{\theta}} \|_{C^{0}}  $$ $$ +\|\partial_s {	u _{+}}	\|_{k-2, p, \delta} \cdot  \|\partial_s u_-\circ \tau_{-2R_\theta} \|_{C^{k-2}}  $$ $$ +\|\partial_s {	u _{+}}	\|_{C^{0}}
\cdot  \|\partial_s u_-\circ \tau_{-2R_{\theta}} \|_{k-1, p, \delta} \}.
$$

$$\sim \|DJ\|_{C^{k-1}}\cdot c(\beta, \gamma) \cdot \| {	u}	\|_{k, p, \delta} \{ \| {	u _{+}}	\|_{k, p, \delta} \cdot  \|\partial_s u_-\circ \tau_{-2R_\theta} \|_{C^{k-2}}  $$ $$ +\|\partial_s {	u _{+}}	\|_{C^{0}}
\cdot  \|\partial_s u_-\circ \tau_{-2R_{\theta}} \|_{k-1, p, \delta} \}.
= \|DJ\|_{C^{k-1}}\cdot c(\beta, \gamma) \cdot \| {	u}	\|_{k, p, \delta}\cdot III.$$

For   $t\in [R-d-l-3 ,R+d+l+3 ]$,
$$ III\sim  \| {	u _{+}}	\|_{k, p, \delta} \cdot  \|\partial_s u_-\circ \tau_{-2R_\theta}|_{[R-d-l-3 ,R+d+l+3 ]} \|_{k-1, p}  $$ $$ +\|\partial_s {	u _{+}}|_{[R-d-l-3 ,R+d+l+3 ]	\|_{1, p}}
\cdot  \|\partial_s u_-\circ \tau_{-2R_{\theta}}|_{[R-d-l-3 ,R+d+l+3 ]} \|_{k-1, p, \delta} $$ $$\leq  \| {	u _{+}}	 \|_{k, p, \delta} \cdot  \|\partial_s u_-|_{[-R-d-l-3 ,-R+d+l+3 ]} \|_{k-1, p}  $$ $$ +\|\partial_s {	u _{+}}|_{[R-d-l-3 ,R+d+l+3 ]	\|_{1, p}}
\cdot e^{\delta (R+d+l+3) } \|\partial_s u_-|_{[-R-d-l-3 ,-R+d+l+3 ]} \|_{k-1, p}  $$

$$ \leq\| {	u _{+}}	\|_{k, p, \delta} \cdot e^{\delta (-R+d+l+3)}   \|\partial_s u_- \|_{k-1, p,\delta} $$ $$ +e^{\delta (-R+d+l+3)}\|\partial_s {	u _{+}}	\|_{1, p, \delta}
\cdot e^{\delta (R+d+l+3) }\cdot e^{\delta (-R+d+l+3)} \|\partial_s u_- \|_{k-1, p, \delta}$$ $$  \sim e^{-\delta R/2}\cdot \| {	u _{-}}	\|_{k, p, \delta} \cdot\| {	u _{+}}	\|_{k, p, \delta}.$$

\QED

By the last lemma in \cite{3}, this finishes the proof of the Theorem 1.1.


\end{document}